\documentclass[letterpaper, 10 pt, conference]{ieeeconf}  

\IEEEoverridecommandlockouts                              
\usepackage{amsfonts}
\usepackage{graphicx,xcolor}
\usepackage{amsmath}
\usepackage{amssymb}
\usepackage{dsfont}

\usepackage{ifthen}
\usepackage{color}
\usepackage{cite}

\usepackage{mathtools}
\usepackage{booktabs}
\usepackage{url}
\usepackage{hyperref}
\usepackage{caption}
\usepackage{subfig}
\usepackage{float}
\usepackage{multirow}

\def\scr#1{{\cal #1}}
\newcommand{\R}{{\rm I\!R}}
\def\eq#1{\begin{equation}#1\end{equation}}

\newcommand{\bbb}{\mathbb}
\newtheorem{theorem}{Theorem}
\newtheorem{lemma}{Lemma}

\newtheorem{assumption}{Assumption}

\newcommand{\dfb}{\stackrel{\Delta}{=}}
\def\qed{ \rule{.08in}{.08in}}
\DeclarePairedDelimiter\ceil{\lceil}{\rceil}
\DeclarePairedDelimiter\floor{\lfloor}{\rfloor}

\newcommand{\1}{\mathbf{1}}
\newcommand{\0}{\mathbf{0}}

\title{\LARGE \bf Heterogeneous Distributed Subgradient
}

\author{Yixuan Lin \hspace{.3in} Ji Liu
\thanks{
Y. Lin is with the Department of Applied Mathematics and Statistics at Stony Brook University (\texttt{yixuan.lin.1@stonybrook.edu}).
J. Liu is with the Department of Electrical and Computer Engineering at Stony Brook University
(\texttt{ji.liu@stonybrook.edu}).
}
}

\begin{document}

\maketitle
\thispagestyle{empty}
\pagestyle{empty}


\begin{abstract}
The paper proposes a heterogeneous push-sum based subgradient algorithm for multi-agent distributed convex optimization in which each agent can arbitrarily switch between subgradient-push and push-subgradient at each time. It is shown that the heterogeneous algorithm converges to an optimal point at an optimal rate over time-varying directed graphs.
\end{abstract}


\section{Introduction}

Stemming from the pioneering work by Nedi\'c and Ozdaglar \cite{nedic2009distributed}, distributed optimization for multi-agent systems has attracted considerable interest and achieved great success in both theory and practice. Surveys of this area can be found in  \cite{yang2019survey,nedic2018distributed,molzahn2017survey}. A typical distributed optimization problem is formulated as follows.

Consider a multi-agent network consisting of $n$ agents, labeled $1$ through $n$ for the purpose of presentation. Every agent is not conscious of such a global labeling, but is capable of distinguishing between its neighbors. The neighbor relations among the $n$ agents are characterized by a possibly time-dependent directed graph $\bbb{G}(t) = (\mathcal{V},\mathcal{E}(t))$ whose
vertices correspond to agents and whose directed edges (or arcs) depict neighbor relations, where $\mathcal{V}=\{1,\ldots,n\}$ is the vertex set and $\mathcal{E}(t)\subset\mathcal{V} \times \mathcal{V}$ is the directed edge set at time $t$.
To be more precise, agent $j$ is an in-neighbor of agent $i$ at time $t$ if $(j,i)\in\scr{E}(t)$, and similarly, agent $k$ is an out-neighbor of agent $i$ at time $t$ if $(i,k)\in\scr{E}(t)$. The directions of arcs represent the directions of information flow in that each agent can send information to its out-neighbors and receive information from its in-neighbors. For convenience, we assume that each agent is always an in- and out-neighbor of itself, implying that $\bbb{G}(t)$ has self-arcs at all vertices for any time $t$. We use $\mathcal{N}_i(t)$ and $\mathcal{N}_i^{-}(t)$ to denote the in- and out-neighbor set of agent $i$ at time $t$, respectively, i.e.,
\begin{align*}
    \mathcal{N}_i(t) &= \{ j \in \mathcal{V}  : ( j, i ) \in \mathcal{E}(t) \}, \\ \mathcal{N}_i^{-}(t) &= \{ k \in \mathcal{V}  :  ( i, k ) \in \mathcal{E}(t) \}.
\end{align*}
It is easy to see that $\mathcal{N}_i(t)$ and $\mathcal{N}_i^{-}(t)$ are always nonempty since they both contain index $i$.
The goal of the $n$ agents is to cooperatively minimize the cost function
$$f(z)=\frac{1}{n}\sum_{i=1}^n f_i(z)$$ in which each $f_i:\R^d\rightarrow\R$ is a ``private'' convex (not necessarily differentiable) function only known to agent $i$. It is assumed that the set of optimal solutions to $f$, denoted by $\scr Z$, is nonempty and bounded.

To solve the distributed optimization problem just described, efforts have been made to design distributed multi-agent versions for various optimization algorithms, including the subgradient method \cite{nedic2009distributed}, alternating direction method of multipliers (ADMM) \cite{wei2012distributed}, Nesterov accelerated gradient method \cite{qu2019accelerated}, and proximal gradient descent \cite{8752033}, to name a few. 
Most existing distributed optimization algorithms require that the underlying communication graph be bi-directional or balanced\footnote{A weighted directed graph is called balanced if the sum of all in-weights equals the sum of all out-weights at each of its vertices \cite{gharesifard2013distributed}.}, which allows a distributed manner to construct a doubly stochastic matrix \cite{metro2,weightbalance}. 
To tackle more general, unbalanced, directed graphs, the push-sum based algorithms have been proposed, with subgradient-push \cite{nedic} being the first one, including notable DEXTRA \cite{xi2017dextra} (a push-sum based variant of the well-known EXTRA algorithm \cite{shi2015extra}) and Push-DIGing \cite{nedic2017achieving}. Another approach to deal with unbalanced directed graphs is called push-pull \cite{pushpull} while its state-of-the-art analysis assumes strongly connectedness at each time instance \cite{timevaryingpushpull}. 
Push-sum is thus the most popular and probably the most powerful existing approach to design distributed (optimization) algorithms over time-varying directed graphs. 


All the existing distributed optimization algorithms are homogeneous in that all the agents in a multi-agent network perform the same (order of) operations. Certain heterogeneity has recently been considered and incorporated in algorithm design. Examples include heterogeneous (uncoordinated) stepsize design for a gradient tracking method \cite{localstepsize}, heterogeneous algorithm picking due to the coexistence of different types of agent dynamics in the network (e.g., a mix of continuous- and discrete-time dynamic agents) \cite{heterogdynamics}, and, particularly popular in machine learning, heterogeneous data training for distributed stochastic optimization \cite{heterogdata}. Notwithstanding this, every agent in these algorithms has to adhere to a single protocol, without theoretical guarantee if any deviation from the protocol occurs. 

With these in mind, this paper aims to design a heterogeneous distributed optimization algorithm in which each agent can change its protocol. To be more precise, the iterative algorithm to be proposed will allow each agent to independently decide its order of operations in any iteration. To illustrate the idea, we focus on the subgradient-push method, and expect that the idea also works for other push-sum based first-order optimization methods.

\section{Subgradient-Push and Push-Subgradient}

We begin with the subgradient-push algorithm proposed in \cite{nedic}. The subgradient method was first proposed in \cite{subgradient} for convex but not differentiable functions. 
For such a convex function $h : \R^d\rightarrow \R$, a vector $g\in\R^d$ is called a subgradient of $h$ at point $x$ if
\eq{
h(y)\ge h(x) + g^\top (y-x) \;\; {\rm for \; all} \;\; y\in\R^d.
\label{eq:subgradient}}
Such a vector $g$ always exists for any $x$ and may not be unique. In the special case when $h$ is differentiable at $x$, the subgradient $g$ is unique and equals 
the gradient of $h$ at $x$. 
From \eqref{eq:subgradient} and the Cauchy-Schwarz inequality, 
\eq{
h(y) - h(x) \ge - G \| y-x\|,
\label{eq:G}}
where $\|\cdot\|$ denotes the 2-norm and $G$ is an upper bound for the 2-norm of the subgradients of $h$ at both $x$ and $y$.

The subgradient-push algorithm is as follows\footnote{The algorithm is written in a different but mathematically equivalent form in \cite{nedic}.}:
\begin{align}
    x_i(t+1) &= \sum_{j\in\scr{N}_i(t)} w_{ij}(t)\Big[x_j(t) - \alpha(t)g_j(t)\Big], \label{eq:pushsub_x}\\  y_i(t+1) &= \sum_{j\in\scr{N}_i(t)} w_{ij}(t)y_j(t),\;\;\;\;\; y_i(0)=1\label{eq:pushsub_y},
\end{align}
where $\alpha(t)$ is the stepsize, $g_j(t)$ is a subgradient of $f_j(z)$ at $x_j(t)/y_j(t)$, and $w_{ij}(t)$, $j\in\scr{N}_i(t)$, are positive weights satisfying the following assumption.  

\vspace{.05in}
 
\begin{assumption}\label{assum:weighted matrix}
There exists a constant $\beta>0$ such that for all $i,j\in\scr V$ and $t$, $w_{ij}(t) \ge \beta$ whenever $j\in\scr{N}_i(t)$. For all $i\in\scr V$ and $t$, $\sum_{j\in\scr{N}_i^{-}(t)} w_{ji}(t) = 1$.
\end{assumption}

\vspace{.05in}

A simple choice of $w_{ij}(t)$ is $1/|\mathcal{N}_j^{-}(t)|$ for all $j\in\scr{N}_i(t)$ which can be easily computed in a distributed manner and satisfies Assumption~\ref{assum:weighted matrix} with $\beta=1/n$. Thus, push-sum based algorithms require that each agent be aware of the number of its out-neighbors.

Let $W(t)$ be the $n\times n$ matrix whose $ij$th entry equals $w_{ij}(t)$ if $j\in\scr{N}_i(t)$ and zero otherwise; in other words, we set $w_{ij}(t)=0$ for all $j\notin\scr{N}_i(t)$.
Assumption~\ref{assum:weighted matrix} implies that $W(t)$ is a column stochastic matrix\footnote{A square nonnegative matrix is called a column stochastic matrix if its column sums all equal one.} 
with positive diagonal entries whose zero-nonzero pattern is compliant with the neighbor graph $\bbb{G}(t)$ for all time $t$. 

In implementation, at each time $t$, each agent $j$ transmits two pieces of information, $w_{ij}(t)[x_j(t) - \alpha(t)g_j(t)]$ and $w_{ij}(t)y_j(t)$, to its out-neighbour $i$, and then each agent $i$ updates its two variables as above. Note that if all $\alpha(t)g_j(t) = 0$, the algorithm simplifies to the push-sum algorithm \cite{push}. Thus, at each time, each agent first performs a subgradient operation, and then follows the push-sum updates. This is why the algorithm \eqref{eq:pushsub_x}--\eqref{eq:pushsub_y} is called subgradient-push. 
It has been recently proved that subgradient-push converges at a rate of $O(1/\sqrt{t})$ over time-varying unbalanced directed graphs, which is the
same as that of the single-agent subgradient and thus optimal \cite{pushcdc}.

Note that in the subgradient-push algorithm, all the agents in a multi-agent network perform the same order of operations, namely an optimization step (subgradient) followed by the push-sum updates. In this paper, we aim to relax this order restriction. To this end, we first introduce a variant of subgradient-push in which the order of subgradient and push-sum operations is swapped. To be more precise, each agent $i$ updates its variables as
\begin{align}
    x_i(t+1) &= \sum_{j\in\scr{N}_i(t)} w_{ij}(t)x_j(t) - \alpha(t)g_i(t), \label{eq:pushfirst_x}\\  y_i(t+1) &= \sum_{j\in\scr{N}_i(t)} w_{ij}(t)y_j(t),\;\;\;\;\; y_i(0)=1\label{eq:pushfirst_y},
\end{align}
where $\alpha(t)$, $w_{ij}(t)$, and $g_i(t)$ are the same as those in subgradient-push.
In the above algorithm \eqref{eq:pushfirst_x}--\eqref{eq:pushfirst_y} each agent~$i$ performs the push-sum updates first for both variables and then the subgradient update for $x_i$ variable. We thus call the algorithm push-subgradient.

Push-subgradient can achieve the same performance as subgradient-push, namely, it converges to an optimal solution at a rate of $O(1/\sqrt{t})$ for general convex functions over time-varying unbalanced directed graphs. It turns out that both push-subgradient and subgradient-push are special cases of the following heterogeneous algorithm.

\section{Heterogeneous Subgradient}

Let $\sigma_i(t)$ be a switching signal of agent $i$ which takes values in $\{0,1\}$. At each time $t$, each agent $j$ transmits two pieces of information, $w_{ij}(t)[x_j(t) - \alpha(t)g_j(t)\sigma_j(t)]$ and $w_{ij}(t)y_j(t)$, to its out-neighbour $i$, and then each agent $i$ updates its variables as follows: 
\begin{align}
    x_i(t+1) =\;& \sum_{j\in\scr{N}_i(t)} w_{ij}(t)\Big[x_j(t) - \alpha(t)g_j(t)\sigma_j(t)\Big] \nonumber\\
    &\;\; - \alpha(t)g_i(t)\big(1-\sigma_i(t)\big),\;\;\; x_i(0)\in\R^d, \label{eq:mix_x}\\  y_i(t+1) =\;& \sum_{j\in\scr{N}_i(t)} w_{ij}(t)y_j(t),\;\;\;\; y_i(0)=1\label{eq:mix_y},
\end{align}
where $\alpha(t)$ is the stepsize, $w_{ij}(t)$, $j\in\scr{N}_i(t)$, are positive weights satisfying Assumption \ref{assum:weighted matrix}. 

In the case when all $\sigma_i(t)=1$, $i\in\scr V$, the above algorithm simplifies to the subgradient-push algorithm \eqref{eq:pushsub_x}--\eqref{eq:pushsub_y}. In the case when all $\sigma_i(t)=0$, $i\in\scr V$, the above algorithm simplifies to the push-subgradient algorithm \eqref{eq:pushfirst_x}--\eqref{eq:pushfirst_y}. 
Thus, the algorithm \eqref{eq:mix_x}--\eqref{eq:mix_y} allows each agent to arbitrarily switch between subgradient-push and push-subgradient at any time, and we hence call it {\em heterogeneous distributed subgradient}.

To state the convergence result of the heterogeneous subgradient algorithm just proposed, we need the following typical assumption and concept. 

\vspace{.05in}

\begin{assumption} \label{assum:step-size}
    The step-size sequence $\{\alpha(t)\}$ is positive, non-increasing, and satisfies $\sum_{t=0}^\infty \alpha(t) = \infty$ and $\sum_{t=0}^\infty \alpha^2(t) < \infty$. 
\end{assumption}

\vspace{.05in}

We say that 
     an infinite directed graph sequence $\{ \bbb{G}(t) \}$ is uniformly strongly connected if there exists a positive integer $L$ such that for any $t\ge 0$, the union graph $\cup_{k=t}^{t+L-1} \bbb{G}(k)$ is strongly connected.\footnote{A directed graph is strongly connected if it has a directed path from any vertex to any other vertex. The union of two directed graphs, $\bbb G_p$ and $\bbb G_q$, with the same vertex set, written $\bbb G_p \cup \bbb G_q$, is meant the directed graph with the same vertex set and edge set being the union of the edge set of $\bbb G_p$ and $\bbb G_q$. Since this union is a commutative and associative binary operation, the definition extends unambiguously to any finite sequence of directed graphs with the same vertex set.}
     If such an integer exists, we sometimes say that $\{ \bbb{G}(t) \}$ is uniformly strongly connected by sub-sequences of length $L$.
It is not hard to prove that the above definition is equivalent to the two popular joint connectivity definitions in consensus literature, namely ``$B$-connected'' \cite{nedic2009distributed_quan} and ``repeatedly jointly strongly connected'' \cite{reachingp1}.

Define $z_i(t)=x_i(t)/y_i(t)$ for all $i\in\scr V$ and $\bar z(t) = \frac{1}{n}\sum_{i=1}^n z_i(t)$. It is easy to see that at the initial time $z_i(0)=x_i(0)$ for all $i\in\scr V$ and $\bar z(0) = \frac{1}{n}\sum_{i=1}^n x_i(0)$.

The following theorem shows that the heterogeneous distributed subgradient algorithm \eqref{eq:mix_x}--\eqref{eq:mix_y} still achieves the optimal rate of convergence to an optimal point. 

\vspace{.05in}

\begin{theorem} \label{thm:bound_everage_n_convex_gene}
    Suppose that $\{ \bbb{G}(t) \}$ is uniformly strongly connected and $\|g_i(t)\|$ is uniformly bounded for all $i$ and~$t$.
\begin{itemize}
    \item[1)] If the stepsize $\alpha(t)$ is time-varying and satisfies Assumption~\ref{assum:step-size}, then with $z^*\in\scr Z$,
\begin{align*}
    \lim_{t\rightarrow\infty}f\bigg(\frac{\sum_{\tau =0}^t \alpha(\tau) \bar z(\tau) }{\sum_{\tau =0}^t \alpha(\tau)}\bigg) &= f(z^*), \\
    \lim_{t\rightarrow\infty}f\bigg(\frac{\sum_{\tau =0}^t \alpha(\tau) z_k(\tau) }{\sum_{\tau =0}^t \alpha(\tau)}\bigg) &= f(z^*), \;\;\;\;\; k\in\scr V.
\end{align*}
    
    \item[2)] If the stepsize is fixed and  $\alpha(t) = 1/\sqrt{T}$ for $T>0$ steps, i.e., $t\in\{0,1,\ldots,T-1\}$, then with $z^*\in\scr Z$,
\begin{align*}
    f\bigg(\frac{\sum_{\tau =0}^{T-1} \bar z(\tau)  }{ T } \bigg) - f(z^*) 
    &\le O\Big(\frac{1}{ \sqrt{T}}\Big), \\
    f\bigg(\frac{\sum_{\tau =0}^{T-1} z_k(\tau)  }{ T } \bigg) - f(z^*) 
    &\le O\Big(\frac{1}{ \sqrt{T}}\Big), \;\;\;\;\; k\in\scr V. 
\end{align*}
\end{itemize}
\end{theorem}

\vspace{.1in}

It is easy to show that the above theorem is a consequence of the following theorem. 

\vspace{.05in}

\begin{theorem} \label{thm:bound_everage_n_convex_bound_gene}
    Suppose that $\{ \bbb{G}(t) \}$ is uniformly strongly connected by sub-sequences of length $L$ and that $\|g_i(t)\|$ is uniformly bounded above by a positive number $G$ for all $i\in\scr V$ and $t\ge 0$.
\begin{itemize}
    \item[1)] 
If the stepsize $\alpha(t)$ is time-varying and satisfies Assumption~\ref{assum:step-size}, then for all $t \ge 0$,
\begin{align}
    & f\bigg(\frac{\sum_{\tau =0}^t \alpha(\tau) \bar z(\tau) }{\sum_{\tau =0}^t \alpha(\tau)}\bigg) - f(z^*) \nonumber \\
    \le\; &  \frac{ \| \bar z(0) - z^* \|^2 + G^2\sum_{\tau =0}^t \alpha^2(\tau)}{2\sum_{\tau =0}^t \alpha(\tau) } 
    \nonumber\\
    &+ \frac{2G\alpha(0) \sum_{i=1}^n \|\bar z(0) - z_i(0) \| }{n\sum_{\tau =0}^t \alpha(\tau) } \nonumber\\
    &+ \frac{32G\sum_{i=1}^n \| x_i(0) \|}{\eta}  \frac{\sum_{\tau =0}^{t-1} \alpha(\tau) \mu^\tau  }{\sum_{\tau =0}^t \alpha(\tau) } \nonumber \\
    &+ \frac{32nG^2}{\eta\mu(1-\mu)}  \frac{\sum_{\tau =0}^{t-1} \alpha(\tau)   (  \alpha(0) \mu^{\frac{\tau}{2}}  +  \alpha(\ceil{\frac{\tau}{2}}) )}{\sum_{\tau =0}^t \alpha(\tau) },\label{eq:bound_timevarying_gene}
\end{align}    
\begin{align}
    & f\bigg(\frac{\sum_{\tau =0}^t \alpha(\tau) z_k(\tau) }{\sum_{\tau =0}^t \alpha(\tau)}\bigg) - f(z^*) \nonumber\\
    \le\; & \frac{ \| \bar z(0) - z^* \|^2 + G^2\sum_{\tau =0}^t \alpha^2(\tau)}{2\sum_{\tau =0}^t \alpha(\tau) } \nonumber \\
    &+ \frac{G\alpha(0) \sum_{i=1}^n (\|\bar z(0) - z_i(0) \| + \|z_k(0)  - z_i(0) \|)}{n\sum_{\tau =0}^t \alpha(\tau) } \nonumber\\
    &+ \frac{32nG^2}{\eta\mu(1-\mu)}  \frac{\sum_{\tau =0}^{t-1} \alpha(\tau)   (  \alpha(0) \mu^{\frac{\tau}{2}}  +  \alpha(\ceil{\frac{\tau}{2}}) )}{\sum_{\tau =0}^t \alpha(\tau) }\nonumber\\
    &+ \frac{32G\sum_{i=1}^n \| x_i(0) \|}{\eta}   \frac{\sum_{\tau =0}^{t-1} \alpha(\tau) \mu^\tau  }{\sum_{\tau =0}^t \alpha(\tau) }  , \;k\in\scr V. \label{eq:bound_timevarying_zi_gene}
\end{align} 
\item[2)] If the stepsize is fixed and  $\alpha(t) = 1/\sqrt{T}$ for $T>0$ steps, then
\begin{align}
    & f\bigg(\frac{\sum_{\tau =0}^{T-1} \bar z(\tau)  }{ T } \bigg) - f(z^*) \nonumber \\ 
    \le \; & \frac{ \| \bar z(0) - z^* \|^2 + G^2 }{ 2\sqrt{T} } + \frac{2G \sum_{i=1}^n \|\bar z(0) - z_i(0) \|}{nT } \nonumber\\
    &+ \frac{32 G \sum_{i=1}^n \| x_i(0)  \| }{\eta(1 - \mu)T}    + \frac{32nG^2}{ \eta \mu (1- \mu)\sqrt{T}}, \label{eq:bound_fixed_gene}\\
    & f\bigg(\frac{\sum_{\tau =0}^{T-1} z_k(\tau)  }{ T } \bigg) - f(z^*) \nonumber\\ 
    \le \; &  \frac{ \| \bar z(0) - z^* \|^2 + G^2 }{ 2\sqrt{T} } 
    + \frac{32 G \sum_{i=1}^n \| x_i(0) \|}{\eta(1 - \mu)T} \nonumber\\
    &+ \frac{G \sum_{i=1}^n (\|\bar z(0)  - z_i(0) \|+\|z_k(0) - z_i(0) \|)}{nT } 
    \nonumber\\
    &+  \frac{32nG^2}{ \eta \mu (1- \mu)\sqrt{T}}, \;\;\; k\in\scr V.  \label{eq:bound_fixed_zi_gene}
\end{align}
\end{itemize}
Here $\eta = \frac{1}{n^{nL}}$ and $\mu=(1-\frac{1}{n^{nL}})^{1/L}$.
\end{theorem}

\vspace{.05in}

Theorem \ref{thm:bound_everage_n_convex_bound_gene} is a generalization of Theorems 2 and 3 in \cite{pushcdc}, so its proof requires a more complicated treatment than those of Theorems 2 and 3 in \cite{pushcdc}. It is not surprising that the bounds given in Theorems 2 and 3 in \cite{pushcdc} are slightly better than those in Theorem \ref{thm:bound_everage_n_convex_bound_gene} here as the former are tailored for a special case. 

Theorem~\ref{thm:bound_everage_n_convex_bound_gene} will be proved in the following subsection.

\subsection{Analysis}

We begin with a property of the $y_i(t)$ dynamics \eqref{eq:mix_y} which is independent of the $x_i(t)$ dynamics \eqref{eq:mix_x}. 
Define a time-dependent $n\times n$ matrix $S(t)$ whose $ij$th entry is \begin{equation}\label{eq:s}
    s_{ij}(t)= \frac{w_{ij}(t)y_j(t)}{y_i(t+1)}=\frac{w_{ij}(t)y_j(t)}{\sum_{k=1}^n w_{ik}(t)y_k(t)}.
\end{equation}
The following lemma guarantees that each $s_{ij}(t)$, and thus $S(t)$, are well defined. 

\vspace{.05in}

\begin{lemma}\label{lemma:y_bound}
    If $\{ \bbb{G}(t) \}$ is uniformly strongly connected, then there exists a constant $\eta>0$ such that $n \ge y_i(t) \ge \eta$ for all $i\in\scr V$ and $t\ge 0$.
\end{lemma}

\vspace{.05in}

The lemma is essentially the same as Corollary~2~(b) in~\cite{nedic}, which further proves that if $\{ \bbb{G}(t) \}$ is uniformly strongly connected by sub-sequences of length $L$, 
$\eta \ge \frac{1}{n^{nL}}$.

It is easy to show that each $S(t)$ is a stochastic matrix\footnote{A square nonnegative matrix is called a row stochastic matrix, or simply stochastic matrix, if its row sums all equal one.}. 
An important property of $S(t)$ matrices is as follows. Let $y(t)$ be a vector in $\R^n$ whose $i$th entry is $y_i(t)$ for all $t\ge 0$.

\vspace{.05in}

\begin{lemma} \label{lemma:rela_y_S}
$ y^\top(t) = y^\top(t+1) S(t)$ for all $t\ge 0$.
\end{lemma}

\vspace{.05in}

{\bf Proof of Lemma~\ref{lemma:rela_y_S}:}
From Assumption \ref{assum:weighted matrix}, $\sum_{i=1}^n w_{ij}(t) = 1$ for any $j\in\scr V$. Then, from \eqref{eq:s},
\begin{align*}
    [y^\top(t+1) S(t)]_j &= \sum_{i=1}^n y_i(t+1) s_{ij}(t) \\
    &= \sum_{i=1}^n y_i(t+1) \frac{w_{ij}(t)y_j(t)}{y_i(t+1)} = y_j(t),
\end{align*}
in which $[\cdot]_j$ denotes the $j$th entry of a column vector. 
\hfill $\qed$

\vspace{.05in}

The above property can be linked to the concept of ``absolute probability sequence'' of the sequence of stochastic matrices $\{S(t)\}$; see Proposition 2 in \cite{pushcdc}. 

To proceed, define 
the following time-dependent quantity: 
\eq{\langle z(t) \rangle \dfb \frac{1}{n} \sum_{i=1}^n y_i(t) z_i(t) = \frac{1}{n}\sum_{i=1}^n x_i(t).\label{eq:<>}} 
Since $y_i(t)>0$ by Lemma~\ref{lemma:y_bound} and $\sum_{i=1}^n y_i(t)=n$, the above quantity is a time-varying convex combination of all $z_i(t)$.   
From update \eqref{eq:mix_x}, for all $i\in \scr V$, 
\begin{align} 
    & z_i(t+1) = \frac{x_i(t+1)}{y_i(t+1)} 
    =  \frac{\sum_{j=1}^n w_{ij}(t)x_j(t) - \alpha(t)g_i(t)}{y_i(t+1)} \nonumber\\
    =\; & \sum_{j=1}^n \frac{ w_{ij}(t)y_j(t) }{y_i(t+1)} z_j(t)-  \frac{\alpha(t)g_i(t)}{y_i(t+1)} \nonumber\\
    =\; & \sum_{j=1}^n s_{ij}(t) z_j(t)-  \frac{\alpha(t)g_i(t)}{y_i(t+1)},\nonumber
\end{align}
which, from Lemma~\ref{lemma:rela_y_S}, leads to
\begin{align} 
    &\langle z(t+1) \rangle 
    =\; \sum_{i=1}^n \frac{y_i(t+1)}{n} z_i(t+1)\nonumber\\
    =\; &\sum_{i=1}^n \frac{y_i(t+1)}{n} \sum_{j=1}^n s_{ij}(t) z_j(t)- \sum_{i=1}^n \frac{y_i(t+1)}{n} \frac{\alpha(t)g_i(t)}{y_i(t+1)} \nonumber\\
    =\; &\sum_{j=1}^n \frac{y_j(t)}{n} z_j(t)- \sum_{i=1}^n \frac{\alpha(t)g_i(t)}{n} \nonumber\\
    =\;& \langle z(t) \rangle -   \frac{\alpha(t)}{n} \sum_{i=1}^n g_i(t).\label{eq:update_<z>_mixed}
\end{align}
It is easy to show that the subgradient-push algorithm \eqref{eq:pushsub_x}--\eqref{eq:pushsub_y} and push-subgradient algorithm \eqref{eq:pushfirst_x}--\eqref{eq:pushfirst_y} share the same $\langle z(t) \rangle$ dynamics as given in \eqref{eq:update_<z>_mixed}. This common dynamics is the basis of the following unified analysis for heterogeneous distributed subgradient.
It is also straightforward to get \eqref{eq:update_<z>_mixed} from equation \eqref{eq:<>}, 
update \eqref{eq:mix_x}, and Assumption \ref{assum:weighted matrix} as follows:
\begin{align*} 
    &\langle z(t+1) \rangle 
    = \frac{1}{n}\sum_{i=1}^n x_i(t+1) 
    \nonumber\\
    =\; & \frac{1}{n}\sum_{i=1}^n \sum_{j=1}^n w_{ij}(t)\Big[x_j(t) - \alpha(t)g_j(t)\sigma_j(t)\Big] \nonumber\\ 
    &- \frac{\alpha(t)}{n}\sum_{i=1}^n g_i(t)\big(1-\sigma_i(t)\big) \nonumber\\
    =\; & \frac{1}{n} \sum_{j=1}^n \Big[x_j(t) - \alpha(t)g_j(t)\sigma_j(t)\Big] - \frac{\alpha(t)}{n}\sum_{i=1}^n g_i(t)\big(1-\sigma_i(t)\big) \nonumber\\
    =\;& \langle z(t) \rangle -   \frac{\alpha(t)}{n} \sum_{i=1}^n g_i(t). 
\end{align*}
The above iterative dynamics of $\langle z\rangle$ can be treated (though not exactly the same) as a single-agent subgradient process for the convex cost function $\frac{1}{n}\sum_{i=1}^n f_i(z)$, which is a critical intermediate step.

The remaining analysis logic is as follows. Using the inequality $\| z_i(t) - z^*\|^2 \le 2 \|\langle z(t) \rangle - z^*\|^2 + 2 \|\langle z(t) \rangle - z_i(t)\|^2$, the analysis is then to bound $\|\langle z(t) \rangle - z^*\|^2$ and $\|\langle z(t) \rangle - z_i(t)\|^2$ separately. For the term $\|\langle z(t) \rangle - z_i(t)\|^2$, since all $z_i$ form a consensus process and $\langle z(t) \rangle$ is always a convex combination of all $z_i(t)$, the term can be bounded using consensus related techniques and relatively easy to deal with. Most analysis will focus on bounding the term $\|\langle z(t) \rangle - z^*\|^2$. It is worth noting that from \eqref{eq:<>},
$\|\langle z(t) \rangle - z^*\|^2=\|\frac{1}{n}\sum_{i=1}^n y_i(t)(z_i(t)-z^*)\|^2=\|\frac{1}{n}\sum_{i=1}^n x_i(t)-z^*\|^2$,
which is the actual Lyapunov function.
Also note that update \eqref{eq:update_<z>_mixed} is equivalent to $\bar x(t+1)=\bar x(t) - \frac{\alpha(t)}{n}\sum_{i=1}^n g_i(t)$ where $\bar x(t) = \frac{1}{n}\sum_{i=1}^n x_i(t)$, which is almost the same as the case of average consensus based subgradient \cite{nedic2009distributed} except that each subgradient $g_i$ is taken at point $z_i$ instead of $x_i$. But this $\bar x$ dynamics is elusive without Lemma~\ref{lemma:rela_y_S}.

To prove Theorem~\ref{thm:bound_everage_n_convex_bound_gene}, we need the following lemmas.

\vspace{.05in}

\begin{lemma} \label{lemma:pushsum_product}
If $\{ \bbb{G}(t) \}$ is uniformly strongly connected,  then for any fixed $\tau\ge 0$, 
    $W(t)\cdots W(\tau+1)W(\tau)$
    will converge to the set 
    $\{v\1^\top : v\in\R^n, \1^\top v=1, v>\0\}$ exponentially fast as $t\rightarrow\infty$.\footnote{We use $\0$ and $\1$ to denote the vectors whose entries all equal to $0$ or $1$, respectively, where the dimensions of the vectors are to be understood from the context. We use $v>\0$ to denote a positive vector, i.e., each entry of $v$ is positive.}
\end{lemma}

\vspace{.05in}

The lemma is essentially the same as Corollary~2~(a) in~\cite{nedic}. If $\{ \bbb{G}(t) \}$ is uniformly strongly connected by sub-sequences of length $L$, Lemma~\ref{lemma:pushsum_product}  implies that there exist constants $c>0$ and $\mu \in [0,1)$ and a sequence of stochastic vectors\footnote{A nonnegative vector is called a stochastic vector if
its entries sum to~$1$.} $ \{ v(t)\}$ such that for all $i,j \in \mathcal{V}$ and $t \ge \tau \ge 0$,
\begin{align}\label{mu}
    \big| \big[W(t)\cdots W(\tau+1)W(\tau)\big]_{ij} - v_i(t) \big|\le c \mu^{t-\tau},
\end{align}
where $[\cdot]_{ij}$ denotes the $ij$th entry of a matrix. 
It has been further shown in \cite{nedic} that $c=4$ and $\mu=(1-\frac{1}{n^{nL}})^{1/L}$.

The following lemma is a generalization of Lemma 8 in \cite{pushcdc}, even though its proof follows the similar flow to that in the proof of Lemma 8 in \cite{pushcdc}.

\vspace{.05in}

 \begin{lemma} \label{lemma:bound_consensus_push_SA_gene} 
    If $\{ \bbb{G}(t) \}$ is uniformly strongly connected by sub-sequences of length $L$ and $\|g_i(t)\|$ is uniformly bounded above by a positive number $G$ for all $i$ and $t$, 
    then for all $t \ge 0$ and $i \in \mathcal{V}$,
\begin{align*}
     &\Big\| z_i(t+1) - \frac{1}{n}\sum_{k=1}^n x_k(t) \Big\| \\  \le\;& \frac{8}{\eta}  \mu^t  \sum_{k=1}^n \| x_k(0) \| + \frac{8 nG}{\eta \mu} \sum_{s=0}^t  \mu^{t-s} \alpha(s).
\end{align*}
If, in addition, Assumption~\ref{assum:step-size} holds, for all $t \ge 0$ and $i \in \mathcal{V}$,
\begin{align*}
     &\Big\| z_i(t+1) - \frac{1}{n}\sum_{k=1}^n x_k(t) \Big\| \\
    \le\;& \frac{8}{\eta}  \mu^t  \sum_{k=1}^n \| x_k(0) \| 
    +  \frac{8nG}{\eta \mu (1-\mu)}   \big(  \alpha(0) \mu^{t/2}  +  \alpha(\ceil{t/2}) \big).
\end{align*}
Here $\eta>0$ and $\mu\in(0,1)$ are constants defined in Lemma~\ref{lemma:y_bound} and \eqref{mu}, respectively.
\end{lemma}

\vspace{.05in}

{\bf Proof of Lemma~\ref{lemma:bound_consensus_push_SA_gene}:}
Define 
$$\epsilon_i(t) \dfb \sum_{j=1}^n w_{ij}(t) g_j(t) \sigma_j(t) + g_i(t) (1-\sigma_i(t))$$ for each $i\in\scr V$ and 
\begin{align*}
 \epsilon(t)\dfb \begin{bmatrix}\epsilon_1^\top(t)\cr\vdots\cr \epsilon_{n}^\top(t)\end{bmatrix}\in\R^{n\times d}.
\end{align*}
Note that 
\begin{align}
    &\sum_{i=1}^n  \| \epsilon_i(t) \|\nonumber\\
    \le\;&  \sum_{i=1}^n \Big(\sum_{j=1}^n w_{ij}(t) \|  g_j(t) \| \sigma_j(t) + \| g_i(t) \|( 1-\sigma_i(t))\Big) \nonumber\\
    \le\;& G \Big(\sum_{i=1}^n \sum_{j=1}^n w_{ij}(t) \sigma_j(t) + \sum_{i=1}^n (1-\sigma_i(t)) \Big) \nonumber\\
    =\;& G \Big( \sum_{j=1}^n \sigma_j(t) + n - \sum_{i=1}^n \sigma_i(t)) \Big) = nG,\label{eq:bound_sum_epsilon}
\end{align}
in which we used the fact that $\sum_{i=1}^n w_{ij}(t)=1$.
Similar to the discrete-time state transition matrix, 
let 
$$\Phi_W(t,\tau)\dfb W(t-1)\cdots W(\tau)$$ 
with $t>\tau$.
From \eqref{eq:mix_x}, 
\begin{align*}
    &x(t+1) 
    = W(t)x(t) -\alpha(t) \epsilon(t) \\
    =\;& \Phi_W(t,0) x(0) - \sum_{l=0}^{t-1} \alpha(l) \Phi_W(t,l+1) \epsilon(l) -\alpha(t) \epsilon(t),
\end{align*}
which implies that
\begin{align}
    &W(t+1) x(t+1)\nonumber\\
    =\;& \Phi_W(t+2,0) x(0) - \sum_{l=0}^{t} \alpha(l) \Phi_W(t+2,l+1) \epsilon(l),\label{eq:h_update_tplus1_gene} \\
    &\1^\top x(t+1) 
     = \1^\top x(0) - \sum_{l=0}^{t} \alpha(l) \1^\top \epsilon(l) \label{eq:h_update_1_gene}.
\end{align}
From Lemma~\ref{lemma:pushsum_product} and \eqref{mu}, there exists a sequence  of stochastic vectors $ \{ \phi(t)\}$ such that for all $i,j \in \mathcal{V}$ and $t\ge s \ge 0$, there holds
$
    | [\Phi_W(t+1,s)]_{ij} - \phi_i(t) |\le 4 \mu^{t-s}
$.
Let $ D(s:t) = \Phi_W(t+1,s) - \phi(t) \1^\top$. 
From \eqref{eq:h_update_tplus1_gene} and \eqref{eq:h_update_1_gene}, 
\begin{align*}
    &W(t+1) x(t+1)  - \phi(t+1)\1^\top x(t+1)\\ 
    =\;& \Phi_W(t+2,0) x(0) - \sum_{l=0}^{t} \alpha(l) \Phi_W(t+2,l+1) \epsilon(l) \\
    &- \phi(t+1) \Big(\1^\top x(0) - \sum_{l=0}^{t}  \alpha(l) \1^\top \epsilon(l)\Big) \\
    =\; & \big(\Phi_W(t+2,0) - \phi(t+1) \1^\top\big) x(0) \\
    &- \sum_{l=0}^{t} \alpha(l)\big( \Phi_W(t+2,l+1) - \phi(t+1) \1^\top\big) \epsilon(l) \\
    =\; & D(0:t+1) x(0) - \sum_{l=0}^{t} \alpha(l) D(l+1:t+1) \epsilon(l),
\end{align*}
which implies that
\begin{align*}
    &x(t+1) = W(t) x(t) - \alpha(t) \epsilon(t) \\
    =\;& \phi(t)\1^\top x(t) + D(0:t) x(0) \\
    &- \sum_{l=0}^{t-1} \alpha(l) D(l+1:t) \epsilon(l) - \alpha(t) \epsilon(t).
\end{align*}
From \eqref{eq:mix_y} and the definition of $D(s:t)$,
$
    y(t+1) = \Phi_W(t+1,0) y(0) = D(0:t) \1 + n \phi(t)
$, or equivalently, $y_i(t+1) = [\Phi_W(t+1,0)\1]_i = [D(0:t) \1]_i + n \phi_i(t)$.
Thus, for all $i\in\scr V$,
\begin{align*}
    & z_i(t+1) - \frac{ x(t)^\top \1}{n} 
    = \frac{x_i(t+1)}{y_i(t+1)} - \frac{x(t)^\top \1}{n}\\
    =\; & \frac{\phi_i(t) x(t)^\top \1 + \sum_{k=1}^n [D(0:t)]_{ik} x_k(0) }{[D(0:t) \1]_i + n \phi_i(t)} - \frac{ x(t)^\top \1}{n}\\
    +\; & \frac{-\sum_{l=0}^{t-1} \alpha(l) \sum_{k=1}^n [D(l+1:t)]_{ik} \epsilon_k(l) - \alpha(t) \epsilon_i(t)}{[D(0:t) \1]_i + n \phi_i(t)}
     \\
    =\; & \frac{n \sum_{k=1}^n [D(0:t)]_{ik} x_k(0)   - [D(0:t) \1]_i x(t)^\top \1 }{n[D(0:t) \1]_i + n^2 \phi_i(t)}\\
    +\; & \frac{-  n\sum_{l=0}^{t-1} \alpha(l) \sum_{k=1}^n [D(l+1:t)]_{ik} \epsilon_k(l)  - n \alpha(t) \epsilon_i(t) }{n[D(0:t) \1]_i + n^2 \phi_i(t)}.
\end{align*}
From Lemma \ref{lemma:y_bound}, $y_i(t+1) \ge \eta$, so is $[D(0:t) \1]_i + n \phi_i(t)$ for all $i\in \scr V$. 
Thus, 
\begin{align*}
    & \Big\| z_i(t+1) - \frac{ x(t)^\top \1}{n} \Big\| \\
    \le\; & \frac{n \| \sum_{k=1}^n [D(0:t)]_{ik} x_k(0) \|  + \| [D(0:t) \1]_i x(t)^\top \1 \| }{n[D(0:t) \1]_i + n^2 \phi_i(t)} \\
    +\; & \frac{   n\sum_{l=0}^{t-1} \alpha(l) \|\sum_{k=1}^n [D(l+1:t)]_{ik} \epsilon_k(l)\|  + n \alpha(t) \|\epsilon_i(t) \|}{n[D(0:t) \1]_i + n^2 \phi_i(t)} \\
    \le\; & \frac{n  (\max_k [D(0:t)]_{ik}) \sum_{k=1}^n  \| x_k(0)\|   }{n[D(0:t) \1]_i + n^2 \phi_i(t)}\\
    + \;& \frac{ n\sum_{l=0}^{t-1} \alpha(l) (\max_k [D(l+1:t)]_{ik}) \sum_{k=1}^n \| \epsilon_k(l)\| }{n[D(0:t) \1]_i + n^2 \phi_i(t)}\\
    + \;& \frac{ \|[D(0:t) \1]_i x(t)^\top \1 \| + n \alpha(t) \| \epsilon_i(t) \| }{n[D(0:t) \1]_i + n^2 \phi_i(t)}\\
    \le\; & \frac{1}{n\eta} \Big( n (\max_k [D(0:t)]_{ik}) \sum_{k=1}^n  \| x_k(0)\|  \\
    +\; & n (\max_k [D(0:t)]_{ik}) \| x(t)^\top \1 \| + n \alpha(t) \| \epsilon_i(t) \| \\
    +\; &  n\sum_{l=0}^{t-1} \alpha(l) (\max_k [D(l+1:t)]_{ik}) \sum_{k=1}^n \| \epsilon_k(l)\|  \| \Big) \\
    \le\; &  \frac{1}{ \eta} \Big[  4 \mu^t  \sum_{k=1}^n \| x_k(0)\| 
    + \sum_{l=0}^{t-1} \alpha(l) 4 \mu^{t-l-1}    \sum_{k=1}^n \| \epsilon_k(l)\| \\
    +\;& \alpha(t) \| \epsilon_i(t) \| + 4 \mu^t \| x(t)^\top \1 \|\Big]. 
\end{align*}
Also, from \eqref{eq:h_update_1_gene}, 
$$
    \| \1^\top x(t+1) \|
     \le \| \1^\top x(0)\| + \| \sum_{l=0}^{t} \alpha(l) \1^\top \epsilon(l)\|. 
$$
Then, from the above inequality,
\begin{align*}
    & \Big\| z_i(t+1) - \frac{ x(t)^\top \1}{n} \Big\| \\
    \le\; & \frac{4}{ \eta} \Big[   \mu^{t}  \sum_{k=1}^n \| x_k(0)\| + \sum_{l=0}^{t-1} \alpha(l)  \mu^{t-l-1} \sum_{k=1}^n  \|\epsilon_k(l)\| \\
    & + \alpha(t) \| \epsilon_i(t) \|  +  \mu^{t} \| \1^\top x(0)\| +  \mu^{t} \Big\| \sum_{l=0}^{t-1} \alpha(l) \1^\top \epsilon(l)\Big\| \Big]\\
    \le\; & \frac{8}{ \eta} \Big[ \mu^{t}  \sum_{k=1}^n \| x_k(0)\| 
    + \sum_{l=0}^{t} \alpha(l)  \mu^{t-l-1} \sum_{k=1}^n \| \epsilon_k(l)\|  \Big]. 
\end{align*}
Using \eqref{eq:bound_sum_epsilon}, it follows that for all $i \in \mathcal{V}$ and $t \ge 0$,
\begin{align*}
    &\Big\| z_i(t+1) - \frac{ x(t)^\top \1}{n} \Big\|  \\
    \le\;& \frac{8}{\eta}  \mu^t  \sum_{k=1}^n \| x_k(0) \| + \frac{8 nG}{\eta \mu} \sum_{s=0}^t  \mu^{t-s} \alpha(s).
\end{align*}
If the stepsize sequence $\{ \alpha(t) \}$ satisfies Assumption~\ref{assum:step-size}, the above inequality further implies that 
\begin{align*}
    & \Big\| z_i(t+1) - \frac{ x(t)^\top \1}{n} \Big\| \\
    \le\; & \frac{8}{\eta}  \mu^t  \sum_{k=1}^n \| x_k(0)  \| 
    +  \frac{8nG}{\eta\mu} \Big(\sum_{s=0}^{\floor{\frac{t}{2}}}  \mu^{t-s} \alpha(s)
    + \sum_{s=\ceil{\frac{t}{2}}}^{t}  \mu^{t-s} \alpha(s) \Big)   \\
    \le\; & \frac{8}{\eta}  \mu^t  \sum_{k=1}^n \| x_k(0) \| 
    +  \frac{8nG}{\eta \mu (1-\mu)}   \big(  \alpha(0) \mu^{t/2}  +  \alpha(\ceil{t/2}) \big).
\end{align*}
This completes the proof.
\hfill$\qed$

\vspace{.05in}

We are now in a position to prove Theorem~\ref{thm:bound_everage_n_convex_bound_gene}.

\vspace{.05in}

{\bf Proof of Theorem~\ref{thm:bound_everage_n_convex_bound_gene}:}
Note that for all $t \ge 0$ and $i,j\in \mathcal{V}$,
\begin{align}
    &\|\langle z(t+1) \rangle - z_i(t+1) \|+\| z_j(t+1) - z_i(t+1) \| \nonumber\\ 
    \le\;& \Big\|\langle z(t+1) \rangle - \frac{1}{n} \sum_{k=1}^n x_k(t)\Big\| +\Big\|z_j(t+1) - \frac{1}{n} \sum_{k=1}^n x_k(t)\Big\| \nonumber \\
    &+ 2\Big\| z_i(t+1) - \frac{1}{n} \sum_{k=1}^n x_k(t)\Big\|\nonumber\\ 
    \le\;& \sum_{j=1}^n \frac{y_j(t+1)}{n} \Big\|  z_j (t+1) - \frac{1}{n} \sum_{k=1}^n x_k(t) \Big\|\nonumber\\
    &+\Big\|z_j(t+1) - \frac{1}{n} \sum_{k=1}^n x_k(t)\Big\|
    + 2\Big\| z_i(t+1) - \frac{1}{n} \sum_{k=1}^n x_k(t)\Big\|\nonumber\\ 
    \le\;& \frac{32}{\eta}  \mu^t  \sum_{i=1}^n \| x_i(0) \| +  \frac{32nG}{\eta \mu}   \sum_{s=0}^t  \mu^{t-s} \alpha(s),\label{eq:bound_for_difference2_general_zi_gene}
\end{align}
where we used Lemma~\ref{lemma:bound_consensus_push_SA_gene} in the last inequality. Similarly, for all $t \ge 0$ and $i\in \mathcal{V}$,
\begin{align}
    &\|\langle z(t+1) \rangle - z_i(t+1) \|+\|\bar z(t+1) - z_i(t+1) \| \nonumber\\ 
    \le\;& \Big\|\langle z(t+1) \rangle - \frac{1}{n} \sum_{k=1}^n x_k(t)\Big\| +\Big\|\bar z(t+1) - \frac{1}{n} \sum_{k=1}^n x_k(t)\Big\| \nonumber\\
    &+ 2\Big\| z_i(t+1) - \frac{1}{n} \sum_{k=1}^n x_k(t)\Big\|\nonumber\\ 
    \le\;& \sum_{j=1}^n \Big(\frac{y_j(t+1)}{n}+\frac{1}{n}\Big) \Big\|  z_j (t+1) - \frac{1}{n} \sum_{k=1}^n x_k(t) \Big\| \nonumber\\
    &+ 2 \Big\| z_i(t+1) - \frac{1}{n} \sum_{k=1}^n x_k(t)  \Big\|\nonumber\\
    \le\;&   \frac{32}{\eta}  \mu^t  \sum_{k=1}^n \| x_k(0) \| +  \frac{32nG}{\eta\mu}   \sum_{s=0}^t  \mu^{t-s} \alpha(s).\label{eq:bound_for_difference2_general_gene}
\end{align}
If, in addition, the stepsize sequence $\{ \alpha(t) \}$ satisfies Assumption~\ref{assum:step-size}, the above two inequalities can be further bounded by Lemma~\ref{lemma:bound_consensus_push_SA_gene} as follows:
\begin{align}
    & \|\langle z(t+1) \rangle - z_i(t+1) \|+\|\bar z(t+1) - z_i(t+1) \|  \nonumber\\
    \le\;& \frac{32}{\eta}  \mu^t  \sum_{k=1}^n \| x_k(0) \|  +  \frac{32nG}{\eta\mu(1-\mu)}   \big(  \alpha(0) \mu^{t/2}  +  \alpha(\ceil{t/2}) \big), \label{eq:bound_for_difference2_gene}\\
    & \|\langle z(t+1) \rangle - z_i(t+1) \|+\|z_j(t+1) - z_i(t+1) \| \nonumber\\
    \le\;& \frac{32}{\eta}  \mu^t \sum_{i=k}^n \|  x_k(0) \| +  \frac{32nG}{\eta\mu(1-\mu)}   \big(  \alpha(0) \mu^{t/2}  +  \alpha(\ceil{t/2}) \big). \label{eq:bound_for_difference2_zi_gene} 
\end{align}
From \eqref{eq:update_<z>_mixed}, for any $z^*\in\scr Z$,
\begin{align}
    & \| \langle z(t+1) \rangle - z^* \|^2
    = \Big\| \langle z(t) \rangle - z^* - \frac{\alpha(t)}{n}\sum_{i=1}^n g_i(t) \Big\|^2 \nonumber\\
    \le\; & \| \langle z(t) \rangle - z^* \|^2 + \Big\| \frac{\alpha(t)}{n}\sum_{i=1}^n g_i(t) \Big\|^2 \nonumber\\
    &- 2 (\langle z(t) \rangle - z^*)^\top \Big(\frac{\alpha(t)}{n}\sum_{i=1}^n g_i(t)\Big) \nonumber \\
    \le\; & \| \langle z(t) \rangle - z^* \|^2 + \alpha^2(t) G^2  \nonumber\\
    &- 2 (\langle z(t) \rangle - z^*)^\top \Big(\frac{\alpha(t)}{n}\sum_{i=1}^n g_i(t)\Big), \label{eq:equation 1_gene}
\end{align}
where we used the convexity of squared 2-norm in the last inequality. 
Moreover, for all $i,k\in\scr V$,
\begin{align}
    & (\langle z(t) \rangle - z^*)^\top g_i(t) \nonumber\\
    =\; & (\langle z(t) \rangle - z_i(t) )^\top g_i(t) + ( z_i(t) - z^*)^\top g_i(t) \nonumber \\
    \ge\; &  f_i(z_i(t))  - f_i(z^*) - G \|\langle z(t) \rangle - z_i(t) \|    \label{eq:equation 4_1_gene} \\
    \ge\; & f_i( z_k(t))  - f_i(z^*) - G \|\langle z(t) \rangle - z_i(t) \|  \nonumber\\
    &- G \| z_k(t) - z_i(t) \|,
    \label{eq:equation 4_zi_gene}
\end{align}
where we used \eqref{eq:subgradient} and \eqref{eq:G} in deriving \eqref{eq:equation 4_1_gene}, and made use of \eqref{eq:G} to get \eqref{eq:equation 4_zi_gene}. 
Similarly, for all $i\in\scr V$,
\begin{align}
     & (\langle z(t) \rangle - z^*)^\top g_i(t) \nonumber\\
    \ge\; & f_i(\bar z(t))  - f_i(z^*) - G \|\langle z(t) \rangle - z_i(t) \|  - G \|\bar z(t) - z_i(t) \|. \label{eq:equation 4_gene}
\end{align}
Combining \eqref{eq:equation 1_gene} and \eqref{eq:equation 4_gene}, 
\begin{align*}
    & \| \langle z(t+1) \rangle - z^* \|^2  \\
    \le\; & \| \langle z(t) \rangle - z^* \|^2 + \alpha^2(t) G^2 - 2\alpha(t) ( f(\bar z(t) )  - f(z^*) )\\
    & + \frac{2G\alpha(t)}{n} \sum_{i=1}^n \big( \|\langle z(t) \rangle - z_i(t) \| + \|\bar z(t) - z_i(t) \| \big),
\end{align*}
which implies that
\begin{align*}
    & 2\alpha(t) ( f(\bar z(t) )  - f(z^*)) \\
    \le \;& \| \langle z(t) \rangle - z^* \|^2 + \alpha^2(t) G^2 - \| \langle z(t+1) \rangle - z^* \|^2 \\
    & + \frac{2G\alpha(t)}{n} \sum_{i=1}^n \big(\|\langle z(t) \rangle - z_i(t) \| + \|\bar z(t) - z_i(t) \| \big).
\end{align*}
Summing this relation over time, it follows that 
\begin{align*}
    & \sum_{\tau =0}^t 2\alpha(\tau) ( f(\bar z(\tau) )  - f(z^*) ) \\
    \le \;& \| \langle z(0) \rangle - z^* \|^2  - \| \langle z(t+1) \rangle - z^* \|^2 + G^2\sum_{\tau =0}^t \alpha^2(\tau) \\
    & + \sum_{\tau =0}^t  \frac{2G\alpha(\tau)}{n} \sum_{i=1}^n \big( \|\langle z(\tau) \rangle - z_i(\tau) \| + \|\bar z(\tau) - z_i(\tau) \| \big).
\end{align*}
Then, 
\begin{align}
    & f\bigg(\frac{\sum_{\tau =0}^t \alpha(\tau) \bar z(\tau) }{\sum_{\tau =0}^t \alpha(\tau)}\bigg) - f(z^*) \nonumber\\
    \le\;& \frac{ \sum_{\tau =0}^t 2\alpha(\tau) ( f(\bar z(\tau)) - f(z^*) ) }{\sum_{\tau =0}^t 2\alpha(\tau)}\nonumber\\
    \le\;& \frac{ \| \langle z(0) \rangle - z^* \|^2  - \| \langle z(t+1) \rangle - z^* \|^2 + G^2\sum_{\tau =0}^t \alpha^2(\tau) }{\sum_{\tau =0}^t 2\alpha(\tau) }  \nonumber\\
    & + \frac{\sum_{\tau =0}^t \frac{2G\alpha(\tau)}{n} \sum_{i=1}^n (\|\langle z(\tau) \rangle - z_i(\tau) \| + \|\bar z(\tau)  - z_i(\tau) \|)}{\sum_{\tau =0}^t 2\alpha(\tau) } \nonumber\\
    \le\;&  \frac{\sum_{\tau =0}^t {G\alpha(\tau)} \sum_{i=1}^n (\|\langle z(\tau) \rangle - z_i(\tau) \| + \|\bar z(\tau)  - z_i(\tau) \|)}{{n}\sum_{\tau =0}^t \alpha(\tau) }  \nonumber\\
    & +  \frac{ \| \langle z(0) \rangle - z^* \|^2 + G^2\sum_{\tau =0}^t \alpha^2(\tau) }{\sum_{\tau =0}^t 2\alpha(\tau) }. \label{eq:bound_mid_gene}
\end{align} 
Similarly, combining \eqref{eq:equation 1_gene} and \eqref{eq:equation 4_zi_gene}, for any $k\in\scr V$, 
\begin{align*}
    & \| \langle z(t+1) \rangle - z^* \|^2 \\
    \le\;& \| \langle z(t) \rangle - z^* \|^2 + \alpha^2(t) G^2 - 2\alpha(t) ( f(z_k(t) )  - f(z^*) )\\
    &  + \frac{2G\alpha(t)}{n} \sum_{i=1}^n \big( \|\langle z(t) \rangle - z_i(t) \| + \| z_k(t) - z_i(t) \| \big),
\end{align*}
which, using the preceding argument, leads to
\begin{align}
    & f\bigg(\frac{\sum_{\tau =0}^t \alpha(\tau) z_k(\tau) }{\sum_{\tau =0}^t \alpha(\tau)}\bigg) - f(z^*) \nonumber\\
    \le\;&  \frac{\sum_{\tau =0}^t G\alpha(\tau) \sum_{i=1}^n (\|\langle z(\tau) \rangle - z_i(\tau) \| + \|z_k(\tau)  - z_i(\tau) \|)}{{n}\sum_{\tau =0}^t \alpha(\tau) }  \nonumber\\
    & +  \frac{ \| \langle z(0) \rangle - z^* \|^2 + G^2\sum_{\tau =0}^t \alpha^2(\tau) }{\sum_{\tau =0}^t 2\alpha(\tau) }. \label{eq:bound_mid_zi_gene}
\end{align} 

We next consider the time-varying and fixed stepsizes separately. 

1) If the stepsize $\alpha(t)$ is time-varying and satisfies Assumption~\ref{assum:step-size}, then combining \eqref{eq:bound_for_difference2_gene} and \eqref{eq:bound_mid_gene}, 
\begin{align*}
    & f\bigg(\frac{\sum_{\tau =0}^t \alpha(\tau) \bar z(\tau) }{\sum_{\tau =0}^t \alpha(\tau)}\bigg) - f(z^*) \\
    \le\;& \frac{ \| \langle z(0) \rangle - z^* \|^2 + G^2\sum_{\tau =0}^t \alpha^2(\tau) }{\sum_{\tau =0}^t 2\alpha(\tau) }  \\
    & + \frac{ G\alpha(0) \sum_{i=1}^n (\|\langle z(0) \rangle - z_i(0) \| + \|\bar z(0)  - z_i(0) \|)}{n \sum_{\tau =0}^t \alpha(\tau) } \\
    & + \frac{32G}{\eta}  \Big(\sum_{i=1}^n \| x_i(0) \|\Big) \frac{\sum_{\tau =0}^{t-1} \alpha(\tau) \mu^\tau  }{\sum_{\tau =0}^t \alpha(\tau) } \\
    & + \frac{32nG^2}{\eta\mu(1-\mu)} \cdot \frac{\sum_{\tau =0}^{t-1} \alpha(\tau)   (  \alpha(0) \mu^{\frac{\tau}{2}}  +  \alpha(\ceil{\frac{\tau}{2}}) )}{\sum_{\tau =0}^t \alpha(\tau) }.
\end{align*} 
Similarly, combining \eqref{eq:bound_for_difference2_zi_gene} and \eqref{eq:bound_mid_zi_gene}, 
\begin{align*}
    & f\bigg(\frac{\sum_{\tau =0}^t \alpha(\tau) z_k(\tau) }{\sum_{\tau =0}^t \alpha(\tau)}\bigg) - f(z^*) \\
    \le\;& \frac{ \| \langle z(0) \rangle - z^* \|^2 + G^2\sum_{\tau =0}^t \alpha^2(\tau) }{\sum_{\tau =0}^t 2\alpha(\tau) }  \\
    & + \frac{ G\alpha(0) \sum_{i=1}^n (\|\langle z(0) \rangle - z_i(0) \| + \|z_k(0)  - z_i(0) \|)}{n\sum_{\tau =0}^t \alpha(\tau) } \\
    & + \frac{32G}{\eta}  \Big(\sum_{i=1}^n \| x_i(0) \|\Big) \frac{\sum_{\tau =0}^{t-1} \alpha(\tau) \mu^\tau  }{\sum_{\tau =0}^t \alpha(\tau) } \\
    & + \frac{32nG^2}{\eta\mu(1-\mu)} \cdot \frac{\sum_{\tau =0}^{t-1} \alpha(\tau)   (  \alpha(0) \mu^{\frac{\tau}{2}}  +  \alpha(\ceil{\frac{\tau}{2}}) )}{\sum_{\tau =0}^t \alpha(\tau) }.
\end{align*} 
Note that $ \langle z(0) \rangle = \frac{1}{n}\sum_i^n z_i(0) = \bar z(0)$. We thus have derived \eqref{eq:bound_timevarying_gene} and \eqref{eq:bound_timevarying_zi_gene}.

2) If the stepsize is fixed and  $\alpha(t) = 1/\sqrt{T}$ for all $t\ge0$, then from \eqref{eq:bound_mid_gene} and \eqref{eq:bound_for_difference2_general_gene},  
\begin{align*}
    & f\bigg(\frac{\sum_{\tau =0}^{T-1} \bar z(\tau)  }{ T } \bigg) - f(z^*) \\
    \le\;& \frac{{G} \sum_{\tau =0}^{T-1} \sum_{i=1}^n \|\langle z(\tau) \rangle - z_i(\tau) \|+\|\bar z(\tau) - z_i(\tau) \|}{{n}T } \\
    & + \frac{ \| \langle z(0) \rangle - z^* \|^2 + G^2 }{ 2\sqrt{T} } \\
    \le\;& \frac{G \sum_{i=1}^n \|\langle z(0) \rangle - z_i(0) \|+\|\bar z(0) - z_i(0) \|}{nT } \\
    & + \frac{ \| \langle z(0) \rangle - z^* \|^2 + G^2 }{ 2\sqrt{T} } + \frac{32 G }{T \eta}  \Big(\sum_{i=1}^n \| x_i(0)\|\Big) \sum_{\tau =0}^{T-2} \mu^\tau\\ 
    &    + \frac{32n G^2 }{T \eta \mu} \sum_{\tau =0}^{T-2} \sum_{s=0}^\tau  \mu^{\tau-s} \frac{1}{\sqrt{T}}\\
    \le\;&  \frac{G \sum_{i=1}^n \|\langle z(0) \rangle - z_i(0) \|+\|\bar z(0) - z_i(0) \|}{{n}T } \\
    & +\frac{ \| \langle z(0) \rangle - z^* \|^2 + G^2 }{ 2\sqrt{T} } \\
    & + \frac{32 G }{T \eta(1 - \mu)} \sum_{i=1}^n \|  x_i(0) \| + \frac{32n G^2 }{ \sqrt{T}\eta \mu (1- \mu)} .
\end{align*}
Similarly, from \eqref{eq:bound_mid_zi_gene} and \eqref{eq:bound_for_difference2_general_zi_gene}, 
\begin{align*}
    & f\bigg(\frac{\sum_{\tau =0}^{T-1} z_k(\tau)  }{ T } \bigg) - f(z^*) \\
    \le\;& \frac{ G \sum_{\tau =0}^{T-1} \sum_{i=1}^n \|\langle z(\tau) \rangle - z_i(\tau) \|+\|z_k(\tau) - z_i(\tau) \|}{nT } \\
    & + \frac{ \| \langle z(0) \rangle - z^* \|^2 + G^2 }{ 2\sqrt{T} } \\
    \le\;& \frac{G \sum_{i=1}^n \|\langle z(0) \rangle - z_i(0) \|+\|z_k(0) - z_i(0) \|}{nT } \\
    & + \frac{ \| \langle z(0) \rangle - z^* \|^2 + G^2 }{ 2\sqrt{T} } + \frac{32 G }{T \eta}  \Big(\sum_{i=1}^n \| x_i(0) \|\Big) \sum_{\tau =0}^{T-2} \mu^\tau\\ 
    &  + \frac{32n G^2 }{T \eta \mu} \sum_{\tau =0}^{T-2} \sum_{s=0}^\tau  \mu^{\tau-s} \frac{1}{\sqrt{T}} \\
    \le\;&  \frac{G \sum_{i=1}^n \|\langle z(0) \rangle - z_i(0) \|+\|z_k(0) - z_i(0) \|}{nT } \\
    & +\frac{ \| \langle z(0) \rangle - z^* \|^2 + G^2 }{ 2\sqrt{T} } \\
    & + \frac{32 G }{T \eta(1 - \mu)}  \sum_{i=1}^n \| x_i(0) \|  
    + \frac{32n G^2 }{ \sqrt{T}\eta \mu (1- \mu)}.
\end{align*}
Since $ \langle z(0) \rangle = \frac{1}{n}\sum_i^n z_i(0) = \bar z(0)$, we have derived \eqref{eq:bound_fixed_gene} and \eqref{eq:bound_fixed_zi_gene}.
\hfill $\qed$

\subsection{A Special Case}

In this subsection, we discuss a special case in which $W(t)$ is a doubly stochastic matrix\footnote{A square nonnegative matrix is called a doubly stochastic matrix if its row sums and column sums all equal one.} 
at all time $t\ge 0$. In this case, it is easy to see from \eqref{eq:mix_y} that $y_i(t)=1$ for all $i\in\scr V$ and $t\ge 0$, and thus $z_i(t)=x_i(t)$ for all $i\in\scr V$ and $t\ge 0$. This observation holds for all push-sum based distributed optimization algorithms studied in this paper as they share the same $y_i(t)$ dynamics which is independent of their $x_i(t)$ dynamics. Then, the subgradient-push, push-subgradient, and heterogeneous subgradient algorithms all simplify to average consensus based subgradient algorithms. Specifically, subgradient-push \eqref{eq:pushsub_x}--\eqref{eq:pushsub_y} simplifies to 
\begin{align}
    x_i(t+1) = \sum_{j\in\scr{N}_i(t)} w_{ij}(t)\Big[x_j(t) - \alpha(t)g_j(x_j(t))\Big], \label{eq:adapt_diffusion}
\end{align}
and push-subgradient \eqref{eq:pushfirst_x}--\eqref{eq:pushfirst_y} simplifies to 
\begin{align}
x_i(t+1) = \sum_{j\in\scr{N}_i(t)} w_{ij}(t)x_j(t) - \alpha(t)g_i(x_i(t)), \label{eq:diffusion_adapt}
\end{align}
which is the ``standard'' average consensus based distributed subgradient proposed in \cite{nedic2009distributed}.
The two updates \eqref{eq:adapt_diffusion} and \eqref{eq:diffusion_adapt} are analogous to the so-called ``adapt-then-combine'' and ``combine-then-adapt'' diffusion strategies in distributed optimization and learning \cite{sayed_magazine}. Thus, in the special case under consideration, the heterogeneous distributed subgradient algorithm \eqref{eq:mix_x}--\eqref{eq:mix_y} simplifies to 
\begin{align*}
    x_i(t+1) =\;& \sum_{j\in\scr{N}_i(t)} w_{ij}(t)\Big[x_j(t) - \alpha(t)g_j(x_j(t))\sigma_j(t)\Big] \\
    &\;\;- \alpha(t)g_i(x_i(t))\big(1-\sigma_i(t)\big),
\end{align*}
which is an average consensus based heterogeneous distributed subgradient algorithm allowing each agent to arbitrarily switch between updates \eqref{eq:adapt_diffusion} and \eqref{eq:diffusion_adapt}.
The preceding discussion implies that the results in this paper apply to the corresponding average consensus based algorithms.

\section{Conclusion}

In this paper, we have proposed a heterogeneous push-sum based subgradient algorithm in which each agent can arbitrarily switch between subgradient-push and push-subgradient, thus subsuming both subgradient-push and push-subgradient as special cases. It has been proved that the proposed  heterogeneous distributed subgradient algorithm converges to an optimal point at an optimal rate over time-varying directed graphs.
The flexibility (arbitrary switching) of the heterogeneous algorithm is expected to be beneficial to protect privacy against an honest-but-curious adversary or an external eavesdropping adversary.
As one future direction, the proposed idea is anticipated to be applicable to other push-sum based optimization algorithms, for example, DEXTRA \cite{xi2017dextra}, Push-DIGing \cite{nedic2017achieving}, and even stochastic gradient push for distributed deep learning \cite{assran19a}. Other future directions include extending the proposed heterogeneous algorithm to cope with more realistic scenarios such as communication delays, asynchronous updating, and package drops.

\bibliographystyle{unsrt}
\bibliography{push}

\end{document}